\documentclass[11pt]{amsart}
\usepackage{graphicx}
\usepackage{amsfonts}
\usepackage{amssymb}
\usepackage{epstopdf}
\newtheorem{theorem}{Theorem}
\newtheorem{corollary}[theorem]{Corollary}
\newtheorem{prop}[theorem]{Proposition}

\title{Bounds on Accumulation Rates of Eigenvalues on Manifolds with Degenerating Metrics}
\author{Jeffrey McGowan}
\begin{document}

\begin{abstract} We consider a family of manifolds with a class of degenerating 
warped product
metrics $g_\epsilon=\rho(\epsilon,t)^{2a}dt^2 +\rho(\epsilon,t)^{2b}ds_M^2$, with 
$M$ compact, $\rho$ homogeneous degree one, $a \le -1$ and $b > 0$.  We study the Laplace operator acting on $L^{2}$ differential $p$-forms and give sharp
accumulation rates for eigenvalues  near the bottom of the essential spectrum
of the limit manifold with metric $g_{0}$.
\end{abstract}

\maketitle
\section{Introduction}
There are many examples of non-compact manifolds which can be thought of as a `limit' of a sequence of compact manifolds.  Particularly nice examples are hyperbolic manifolds in dimensions 2 and 3; the cusp closing theorem of Thurston \cite{Thurston} then says that every complete, non-compact manifold $M_0$ is the limit of a sequence of hyperbolic manifolds $M_k \to M_0$.  Since the Laplacian on $M_0$ has continuous spectrum, one expects the eigenvalues of $M_k$ to accumulate.  In dimension 2, Ji, Zworski, and Wolpert (\cite{Ji, JiZworski, Wolpert1, Wolpert2}) have given bounds for the accumulation rate of eigenvalues near the bottom of the essential spectrum in the hyperbolic case, while in dimension 3 analogous results were obtained by Chavel and Dodziuk (\cite{ChavelDodziuk}).  Dodziuk and McGowan obtained similar results for the Laplacian acting on differential forms (\cite{DM}).

Colbois and Courtois considered convergence of eigenvalues below the bottom of the essential spectrum in a much more general setting \cite{CC}.  The accumulation rate for eigenvalues of the Laplacian on functions for manifolds $N=\tilde{N}\cup (M^n\times I)$ with 'pseudo-hyperbolic' metrics on $(M^n\times I)$ was given by Judge \cite{Judge}.  Judge also computes the essential spectrum for a more general class of degenerating metrics, and investigates the convergence of eigenfunctions.

We will consider manifolds $N_{\epsilon}=\tilde{N}\cup (M^n\times I)$, $\tilde{N}$ and $M^n$ compact, with $n=dim(M)$, and a family of metrics \begin{equation}\label{metric}g_\epsilon=\rho(\epsilon,t)^{2a}dt^2 +\rho(\epsilon,t)^{2b}ds_M^2\end{equation} on $M^n \times I$.  Here $\rho = c_1\epsilon +c_2t$, $c_1,c_2>0$, $t\in I=[0,1]$, $a\leq -1$, $b>0$, and $ds_M^2$ is the metric on $M^n$.  We identify the boundary of $\tilde{N}$ with $M^{n}\times {1}$.  These are the metrics discussed by Melrose in \cite{Melrose} and considered by Judge in \cite{Judge}.  We consider only non-negative values of $t$ with $t \in [0,1]$, which simplifies the statements of the results, although we must consider manifolds with boundary. 
 The condition $a \leq -1$ means that the limiting manifold $N_{0}$ is complete.

We study the accumulation rate for eigenvalues near the bottom of the essential spectrum of the Laplacian acting on both functions and differential forms.  Our main results are

\begin{theorem}\label{th2} Suppose $N_\epsilon=\tilde{N}\cup (M^n\times I)$, $\tilde{N}$ and $M^n$ compact, with metric \begin{equation}g_\epsilon=\rho(\epsilon,t)^{2a}dt^2 +\rho(\epsilon,t)^{2b}ds_M^2\end{equation} on $M^n \times I$, with $\rho$ as above.  Let $$R=\int_0^1\rho(\epsilon,s)^a\,ds$$ be the geodesic distance from the boundary ${0} \times M^{n})$ of $N_{\epsilon}$ to $\tilde{N}$. Let $\Xi_\epsilon(x^2)$ be the number of eigenvalues of the Laplacian acting on coexact $p$-forms (satisfying absolute boundary conditions on the boundary of $N_{\epsilon}$) in $[\sigma,\sigma + x^2)$ where $\sigma$ is the bottom of the essential spectrum for coexact forms of degree $p$ and $0 < p < n$.  Then $$\Xi_\epsilon(x^2) =\frac{dxR}{\pi} + O_x(1)$$ where $d$ is the dimension of the space of harmonic forms of degree $p$ on $M$.\end{theorem}
This agrees with the results of Judge \cite{Judge}, Chavel and Dodziuk \cite{ChavelDodziuk} and Dodziuk and McGowan \cite{DM} in the special cases they considered.
\begin{theorem}\label{th1} Suppose $N_\epsilon$ is as in Theorem \ref{th1}.  Then the essential spectrum of the Laplacian acting on coexact $p$-forms, $0 \leq p \leq n$ on $N_0$ is $$\begin{array}{lcr}
\left[\left(\frac{n-2p}{2}\right)^2c_2^2b^2,\infty\right)&\qquad&a=-1\\\\
\left[0,\infty\right)&\qquad&a<-1
\end{array}$$\end{theorem}
Note that this agrees with Judges results (\cite{Judge}) for functions when $p=0$, and with Mazzeo and Phillips results for the essential spectrum on geometrically finite hyperbolic manifolds (\cite{MazzeoPhillips}, with $c_2=b=1$ and $a=-1$).  We have recently learned that these results for the essential spectrum have been obtained independently by Antoci (\ref{Antoci}).

This paper is organized as follows.  In Section \ref{geom} we discuss the geometry of the manifolds under consideration, and rewrite the metric (\ref{metric}) in a way which makes the geometry more evident.  In Section \ref{functions} we illustrate our techniques by computing the essential spectrum and accumulation rates for eigenvalues as $\epsilon \to 0$ in the case of functions ($p=0$).  In Section \ref{upperforms} we compute the essential spectrum and give lower bounds on the accumulation rate in the $p\neq 0$ case.  Finally, in Section \ref{lowerforms} we give upper bounds on the accumulation rate for $p \neq 0$, completing the proof of Theorem \ref{th1}.  

We wish to thank J\'ozef Dodziuk for many helpful conversations.

\section{The Geometry}\label{geom}

Metrics of the type (\ref{metric}) are discussed by Melrose in \cite{Melrose}.  When $a \leq -1$ such metrics are complete on the limit manifold $N_{0}$.  Melrose classifies metrics where $a=-1$, $b=1$ as 'hc', or hyperbolic cusp metrics, and metrics where $a=-1$, $b=0$ as 'boundary', or metrics with cylindrical end.  Since we will consider metrics where $a \leq -1$, $b > 0$ we rewrite the metric to make the geometry more evident.

Let $\tau$ be the geodesic distance from $t=0$ to $t=1$, in other words the geodesic distance from the a point $(0,p \in M)$ to $\tilde{N}$.  Then \begin{equation}\label{taudef}
\tau = \int_0^t{\rho(\epsilon,s)^a\,ds}= \int_0^t{(c_1\epsilon + c_2s)^a\,ds}
\end{equation}
and we have two distinct cases, $$\begin{array}{lcc}
\tau =\frac{1}{c_2}\left( \ln\left(\frac{{c_1\epsilon +c_2t}}{c_1\epsilon}\right)\right)&\qquad&a=-1\\
\tau = \frac{1}{c_2}\left(\frac{(c_1\epsilon+c_2t)^{a+1}-(c_1\epsilon)^{a+1}}{c_2(a+1)}\right)&\qquad&a<-1
\end{array}$$

Solving for $t$ and substituting into the metric (\ref{metric}) we get 
\begin{equation}\label{mymetrics}
\begin{array}{lcc}
ds^2=d\tau^2+(c_1\epsilon)^{2b}e^{2bc_2\tau}ds_M^2&\qquad&a=-1\\
ds^2=d\tau^2+(c_2(a+1)\tau+(c_1\epsilon)^{a+1})^\frac{2b}{a+1}ds_M^2&\qquad&a<-1
\end{array}
\end{equation}
which is of the form $ds^2=d\tau^2+f_\epsilon(\tau)ds_M^2$ in both cases.  As $\epsilon \to 0$, $\tau\to \infty$, and we have a warped product $I \times_{f_\epsilon} M$, with the length of the interval given by \begin{equation}\label{rcalc}
\tau(1)=\Bigg\{\begin{array}{lr}R = \frac{1}{c_2}\ln\left(\frac{c_1\epsilon+c_2}{c_1\epsilon}\right)&a=-1\\R = \frac{(c_1\epsilon+c_2)^{a+1}-(c_1\epsilon)^{a+1}}{c_2(a+1)}&a<-1\end{array}
\end{equation}

When $a=-1$, $f_e(\tau)$ gives an essentially hyperbolic metric; one thinks of pinching off a closed geodesic, with $\epsilon$ the length of that geodesic.  When $a < -1$, the cross sections $M$ shrink at a slower rate as one recedes from $\tilde{N}$, and the warped product $I \times_{f_\epsilon} M$ is intermediate between a hyperbolic cusp and a flat cylinder.

Clearly, for any fixed $\epsilon$ and $\tau$, the cross section $\{\tau\}\times f_\epsilon(\tau)M$  has injectivity radius bounded below by some constant.  Moreover, $f_\epsilon(\tau)$ is an increasing function of $\tau$ for all $a \le -1$.  Since $M$ is compact,a scaling argument shows that the first non-zero eigenvalue of the Laplacian acting on coexact forms of degree $p$ on $M_{\epsilon,\tau}$, say $\nu_{p,\epsilon}(\tau)$, is a decreasing function of $\tau$.  Hence, as the geodesic distance of a given cross section from $\tilde{N}$ increases, $\nu_{p,\epsilon}(\tau)$ increases.  This allows us, for technical reasons, to restructure our decomposition of $N$ as follows;$$N=\tilde{N}'\cup\left(M\times [0,R-r_0+1]\right),$$ where $\tilde{N}' = \tilde{N} \cup M \times [R-r_0,R]$.  $r_0$ will be chosen so that $\nu_{p,\epsilon}(r_0)$ is relatively large. 

\section{Functions}\label{functions}

We follow essentially the argument in \cite{DM}.  First, we choose a function $f$ whose restriction to $\tilde{N}'$ is orthogonal to a basis of eigenfunctions with eigenvalues less than or equal to $\sigma + x^2$.  This will only change the counting function $N_\epsilon(x^2)$ by a bounded amount which can be absorbed into the $O_x(1)$ term (\cite[Lemma 3.6]{ChavelDodziuk}).  Next, we decompose $f$ on $M\times [0,R-r_0]$ as $f=\bar{f}+\bar{\bar{f}}$, where $\bar{f}$ depends only on $\tau$ and $\bar{\bar{f}}$ is orthogonal to constants on $M$.  $\bar{f}$ is computed by averaging over each cross section.

Now, if we choose $r_0$ so that $\nu_{p,\epsilon}(r_0)>\sigma + x^2$, then $\bar{\bar{f}}$ does not contribute to the counting function $N_\epsilon(x^2)$.  Concentrating on $\bar{f}$, a straightforward calculation shows that \begin{equation}
\label{funSL}
\Delta\bar{f}=*d*d\bar{f}=\frac{1}{f_\epsilon^{\frac{n}{2}(\tau)}}\frac{d}{d\tau}\left(\frac{d\bar{f}}{d\tau}f_\epsilon^\frac{n}{2}(\tau)\right)
\end{equation}
 is a classical Sturm-Liouville problem, and we can convert to the form (see \cite{CourantHilbert1}) $$u''-ru=\lambda u$$ with \begin{eqnarray*}
u& = & f_\epsilon^\frac{n}{4}(\tau)\bar{f} \\
r & = &\frac{( f_\epsilon^\frac{n}{4}(\tau))''}{f_\epsilon^\frac{n}{4}(\tau)}
\end{eqnarray*}
When $a=-1$, $$r=\left(\frac{nbc_2}{2}\right)^2,$$ and (\ref{funSL}) becomes $$u''=\left(\lambda + \left(\frac{nbc_2}{2}\right)^2\right)u.$$  We get
\begin{prop}\label{fun1}
Suppose $N_\epsilon=\tilde{N}\cup (M^n\times I)$, $\tilde{N}$ and $M^n$ compact, with metric \begin{equation}g_\epsilon=\rho(\epsilon,t)^{-2}dt^2 +\rho(\epsilon,t)^{2b}ds_M^2\end{equation} on $M^n \times I$, with $\rho = c_1\epsilon + c_2 t$.  Then the essential spectrum of the Laplacian acting on functions on $N_0$ is $$
\left[\left(\frac{nc_2b}{2}\right)^2=\sigma,\infty\right).$$ Let  $R$ be as in (\ref{rcalc}) and let $N_\epsilon(x^2)$ be the number of eigenvalues of the Laplacian acting on function  in $[\sigma,\sigma + x^2)$.  Then $$N_\epsilon(x^2) = \frac{xR}{\pi}+O_x(1).$$ \end{prop}
This is as in \cite{Judge}, with slightly different notation.

When $a<-1$, \begin{equation}r=\frac{(a+1)bn(bn-2a-2)c_{2}^{2}}{2(c_{2}(a_+1)\tau+(c_{1}\epsilon)^{a+1})^{2}}.\label{integrable}\end{equation}  The potential (\ref{integrable}) is integrable, and \cite[Theorem 4.1]{ChavelDodziuk} tells us the counting function for the corresponding Sturm-Liouville problem has the same asymptotics as if the potential were identically 0.  Hence,
\begin{prop}
Suppose $N$ is as in Proposition \ref{fun1} with metric \begin{equation}g_\epsilon=\rho(\epsilon,t)^{2a}dt^2 +\rho(\epsilon,t)^{2b}ds_M^2\end{equation} on $M^n \times I$, where $a < -1$ and $\rho$ as above.  Then the essential spectrum of the Laplacian acting on functions on $N_0$ is $
\left[0,\infty\right).$ Let $R$ be as in (\ref{rcalc}) and let $N_\epsilon(x^2)$ be the number of eigenvalues of the Laplacian acting on function  in $[0,x^2)$.  Then \begin{equation}\label{Nfun}N_\epsilon(x^2) = \frac{xR}{\pi}+O_x(1).\end{equation} \end{prop}
The essential spectrum in this case was given by Judge (\cite{Judge}).  The accumulation estimate (\ref{Nfun}) can also be obtained using Judge's techniques (\cite{Judge2}).

\section{Upper eigenvalue bounds for forms}\label{upperforms}

We consider the sequence of eigenvalues of the Laplacian acting on coexact forms of degree $p$, $0 < \nu_{1} \le \nu_{2}\le \cdots \to \infty$.  If we can give an upper bound  $y\ge \nu_{j}$ for some $j$, we will obtain a lower bound for the counting function $\Xi_{\epsilon}(y)  \le j$.  We work in the space $\mathcal{E}$ of $C^{\infty}$ coexact forms of degree $p$ on $N_{\epsilon}$ with support contained in $\{x|1 \le d(x,\tilde{N})\le R\}$, with coefficients which depend only on $\tau$.   Any form $\omega \in \mathcal{E}$ is zero on $\tilde{N}$, and we choose forms \begin{equation}\label{omega}\omega=\sum_{i=1}^{d}b_{i}d\tau\wedge H_{i}\end{equation} where $d$ is the dimension of the space of harmonic $p$ forms on the cross section $M$ and $H_{i},i=1,2,\ldots,d$ is a basis of harmonic $p$ forms on $M$.

Using Courant's min-max principle the eigenvalues $\nu_{j}$ are no greater than the critical values of the Rayleigh-Ritz quotient $(\Delta \omega,\omega)/(\omega,\omega))$ with $\omega \in \mathcal{E}$.  Since $\omega$ is coexact, we have $$\frac{(\Delta \omega,\omega)}{(\omega,\omega))}=\frac{(d\omega,d\omega)}{(\omega,\omega)}$$  Since the $b_{i}$ depend only on $\tau$, and an application of $d$ to the sum in (\ref{omega}) involves only derivatives with respect to other basis elements, we compute $$d{\omega} = \sum_{i=1}^{d}b_{i}'d\tau\wedge H_{i}$$ where the prime indicates differentiation with respect to $\tau$.  Computing the respective $L^{2}$ norms we have \begin{eqnarray}
(\omega,\omega)&=&C_M\sum_{i=1}^{d}\int_{1}^{R}b_{i}^{2}f_{\epsilon}^{\frac{n-2p}{2}}(\tau)\,d\tau\label{bottom}\\
(d\omega,d\omega)&=&C_M\sum_{i=1}^{d}\int_{1}^{R}(b_{i}')^{2}f_{\epsilon}^{\frac{n-2p}{2}}(\tau)\,d\tau
\label{top}\end{eqnarray}  $C_{M}$ can be computed by integrating the basis elements of the cross section $M$.

Using integration by parts in the numerator of the Rayleigh-Ritz quotient, we get $d$ copies of a Sturm-Liuoville problem very similar to the one in section \ref{functions},$$-\frac{1}{f_{\epsilon}^{\frac{n-2p}{2}}}\left(b_{i}'f_{\epsilon}^{\frac{n-2p}{2}}\right)'=\lambda b_{i}.$$  As in Section \ref{functions}, we reduce to the form $$u''-ru=\lambda u.$$  In the pseudo-hyperbolic case, when $a=-1$, we have $$r=\left(\frac{n-2p}{2}\right)^{2}c_{2}^{2}b^{2}$$ and we get \begin{prop}\label{lowerform} Suppose $N_\epsilon=\tilde{N}\cup (M^n\times I)$, $\tilde{N}$ and $M^n$ compact, with metric \begin{equation}g_\epsilon=\rho(\epsilon,t)^{-2}dt^2 +\rho(\epsilon,t)^{2b}ds_M^2\end{equation} on $M^n \times I$, with $\rho=c_{1}\epsilon+c_{2}t$.  If $N_\epsilon(x^2)$ is the number of eigenvalues of the Laplacian acting on coexact $p$-forms in $[\left(\frac{n-2p}{2}\right)^{2}c_{2}^{2}b^{2},\left(\frac{n-2p}{2}\right)^{2}c_{2}^{2}b^{2}+ x^2)$ , then $$N_\epsilon(x^2) \ge \frac{dxR}{\pi} + O_x(1)$$ where $d$ is the dimension of the space of harmonic forms of degree $p$ on $M$.  In this case, $R = \frac{1}{c_2}\ln\left(\frac{c_1\epsilon+c_2}{c_1\epsilon}\right)$, and letting $\epsilon \to 0$, we see that the essential spectrum  of the Laplacian acting on coexact $p$-forms, $0 \leq p \leq n$ on $N_0$ is $\left[\left(\frac{n-2p}{2}\right)^2c_2^2b^2,\infty\right)$ if $d \ne 0$.\end{prop}  

When $a < -1$ a messy but straightforward calculation give $$r=\frac{c_{2}^{2}b(n-2p)(\frac{b(n-2p)}{2}-1}{2(c_{2}(a+1)\tau + (c_{1}\epsilon)^{a+1})^{2}}.$$  This is an integrable potential, and we get \begin{prop}\label{lowerform2} Suppose $N_\epsilon=\tilde{N}\cup (M^n\times I)$, $\tilde{N}$ and $M^n$ compact, with metric \begin{equation}g_\epsilon=\rho(\epsilon,t)^{2a}dt^2 +\rho(\epsilon,t)^{2b}ds_M^2\end{equation} on $M^n \times I$, with $\rho=c_{1}\epsilon+c_{2}t$.  If $N_\epsilon(x^2)$ is the number of eigenvalues of the Laplacian acting on coexact $p$-forms in $[0, x^2)$ , then $$N_\epsilon(x^2) \ge \frac{dxR}{\pi} + O_x(1)$$ where $d$ is the dimension of the space of harmonic forms of degree $p$ on $M$.  In this case, $R = \frac{(c_1\epsilon+c_2)^{a+1}-(c_1\epsilon)^{a+1}}{c_2(a+1)}$, and letting $\epsilon \to 0$, we see that the essential spectrum  of the Laplacian acting on coexact $p$-forms, $0 \leq p \leq n$ on $N_0$ is $\left[0,\infty\right)$ if $d \ne 0$.\end{prop}  

\section{Lower eigenvalue bounds for forms}\label{lowerforms}

We will use the method of \cite[Lemma ?]{McGowan} to get global lower eigenvalue bounds for forms on $N_{\epsilon}$ based on lower eigenvalue bounds on local eigenvalue bounds on (overlapping) pieces of $N_{\epsilon}$.  In particular, we use the idea of constructing a globally defined form while keeping control of the Rayleigh-Ritz quotient as in \cite{DM}.  For details on the underlying \v{C}ech-de\ Rham formalism see \cite[Chapter 2]{BottTu}.  The pieces we will consider might have mildly singular boundaries, but all the familiar results of Hodge theory hold (\cite{Cheeger1,Cheeger2}, see also \cite[Section ?]{McGowan} and \cite[Section 4]{DM}).  We omit many details here, but refer the reader especially to \cite{DM} if they wish to fill in the blanks.

First, we pick a simple open cover of $N_{\epsilon}$ consisting of two pieces; $U_{1}=\tilde{N}'\setminus\partial\tilde{N}'$, and $U_{2}=M\times [0,R-r_{0}+1]$.  Recall that $\tilde{N}' = \tilde{N} \cup M \times [R-r_0,R]$, so $U_{1}$ and $U_{2}$ overlap, with $U_{1}\cap U_{2}=M\times [R-r_{0},R-r_{0}+1]$.  Next, we choose a coexact $p$ form $\phi$ so that the restriction $\phi|_{U_{1}}=\phi_{1}$ is orthogonal to the finite dimensional space of exact eigenforms (on $U_{1}$) of degree $p+1$ with eigenvalue less than or equal to $y^{2}$.  This is possible using \cite[Proposition 5.1]{DM}; the proof must be modified somewhat to account for the more general setting here, but the modifications are simple if messy.  We will specify values for $r_{0}$ and $y^{2}$ later.

Now, since $\phi_{1}$ is assumed to be exact with eigenvalue greater than or equal to $y^{2}$, there exists a unique coexact form $\psi_{1}$ of degree $p$ on $U_{1}$ with $d\psi_{1}=\phi_{1}$ and $$\frac{(\phi_{1},\phi_{1})}{(\psi_{1},\psi_{1})}=\frac{(d\psi_{1},d\psi_{1})}{(\psi_{1},\psi_{1})}\ge y^{2}.$$  Likewise, by exactness, there exists a unique coexact form $\psi_{2}$ on $U_{2}$ with $d\psi_{2}=\phi_{2}$, but we do not yet have any bounds on the Rayleigh-Ritz quotient (and hence on eigenvalues) \begin{equation}\label{RR}\frac{(\phi_{2},\phi_{2})}{(\psi_{2},\psi_{2})}=\frac{(d\psi_{2},d\psi_{2})}{(\psi_{2},\psi_{2})}.\end{equation}

Next, we wish to decompose $\phi_{i}$, $i=1,2$ on $M\times [0,R]$ in a similar fashion to our decomposition for functions at the beginning of Section \ref{functions}.  We will model ourselves on the argument in \cite{DM}, but since the cross section ${\tau} \times M^{n}$ is arbitrary here, we cannot just average coefficients.  Rather, we use a harmonic projection.  First, decompose $\phi_{i} = \alpha \wedge d\tau + \beta$, where $\beta$ does not contain $d\tau$.  Next, use harmonic projection on $\alpha$ to get, $$\alpha = \sum_{i=1}^{d}a_{i}d\tau\wedge H_{i}  + \gamma$$ where $d$ is the dimension of the space of harmonic $p$ forms on $M$, $H_{i}$ is a basis of harmonic $p$ forms on $M$, and the $a_{i}$ depend only on $\tau$.

Now, we can write $\phi_{i}=\bar{\phi_{i}}+\bar{\bar{\phi_{i}}}$, with \begin{eqnarray}
\label{decomp}
\bar{\phi_{i}} &=&  \sum_{i=1}^{d}a_{i}d\tau\wedge H_{i}\\
\bar{\bar{\phi_{i}}} &=& \phi_{i}-\bar{\phi_{i}}  = \beta+\gamma
\end{eqnarray}
We do the same for $\psi_{i}$.  By construction, the coefficients of $\bar{\phi_{i}}$ and $\bar{\psi_{i}}$ depend only on $\tau$.  A straightforward calculation (see, for example, \cite{Dodziuk}), shows that as the metric on the cross sections scales by a factor $f_{\epsilon}(\tau)$, the Rayleigh-Ritz quotient scales as $$\frac{(d\bar{\bar{\phi_{i}}},d\bar{\bar{\phi_{i}}})|_{g_{_{\epsilon}}}} {(\bar{\bar{\phi_{i}}},\bar{\bar{\phi_{i}}})|_{g_{_{\epsilon}}}} = \frac{1}{f_{\epsilon}(\tau)}\left(\frac{(d\bar{\bar{\phi_{i}}},d\bar{\bar{\phi_{i}}})|_{g_{_{1}}}}{(\bar{\bar{\phi_{i}}},\bar{\bar{\phi_{i}}})|_{g_{_{1}}}}\right).$$  For small $\tau$, $f_{\epsilon}(\tau)$ is small, and thus if $r_{0}$ is chosen appropriately, $\bar{\bar{\phi_{i}}}$ will not contribute to any accumulation of eigenvalues. 

So far, we have put only a finite number of conditions, depending only on $x$, on our original selection of $\phi$.  These conditions guarantee that $\phi_{1}$ is orthogonal to the finite dimensional space of exact eigenforms (on $U_{1}$) of degree $p+1$ with eigenvalue less than or equal to $y^{2}$.   We  still need to determine how many additional choices we must make to gain control of the Rayleigh-Ritz quotient (\ref{RR}).   By construction, we can write $$\bar{\phi_{2}}=\sum_{i=1}^{d}a_{i}d\tau\wedge H_{i}$$ where $d$ is the dimension of the space of harmonic $p$ forms on $M$, $H_{i}$ is a basis of harmonic $p$ forms on $M$, and the $a_{i}$ depend only on $\tau$.  Consequently, we can write $$\bar{\psi_{2}}=\sum_{i=1}^{\zeta}f_{i}\,d\tau\wedge\alpha_{i,p-1}+\sum_{i=1}^{d}b_{i}H_{i}$$ with $d\bar{\psi_{2}}=\bar{\phi_{2}}$, $a_{i}=b_{i}'$, and the prime denoting differentiation with respect to $\tau$.

To evaluate the Rayleigh-Ritz quotient on $U_{2}$, we use   \begin{eqnarray}
(\bar{\psi_{2}},\bar{\psi_{2}})&=&C_M\sum_{i=1}^{d}\int_{0}^{R-r_{0}}b_{i}^{2}f_{\epsilon}^{\frac{n-2p}{2}}(\tau)\,d\tau\\
(\bar{\phi_{2}},\bar{\phi_{2}})&=&C_M\sum_{i=1}^{d}\int_{0}^{R-r_{0}}(b_{i}')^{2}f_{\epsilon}^{\frac{n-2p}{2}}(\tau)\,d\tau
\end{eqnarray} and we again have $d$ copies of $$-\frac{1}{f_{\epsilon}^{\frac{n-2p}{2}}}\left(b_{i}'f_{\epsilon}^{\frac{n-2p}{2}}\right)'=\lambda b_{i}.$$   Letting $\sigma$ be the bottom of the essential spectrum for $N_{0}$ we see that the number of eigenvalues in the interval $[\sigma,\sigma+x^{2})$ for the equation $\Delta_{p}\psi_{2}= \nu \psi_{2}$ is given by $\frac{dxR}{\pi}+O_{x}(1)$.  Thus, we can choose $\phi$ in such a way that $\psi_{2}$ is orthogonal in $L^{2}$ to the basis of eigenforms with eigenvalues less than $\sigma + x^{2}$ on $U_{2}$ by imposing $\frac{dxR}{\pi}+O_{x}(1)$ conditions.  The number of conditions imposed on the choice of $\phi$ depends only on $x$.

We have chosen $\phi$ in such a way that we have control over the relevant Rayleigh-Ritz quotients on both $U_{1}$ and $U_{2}$, but it is not the case that $\psi_{1}$ and $\psi_{2}$ must match on $U_{1}\cap U_{2}$.  Since $d\psi_{1}=\phi|_{U_{1}}$ and $d\psi_{2}=\phi|_{U_{2}}$ it is clear that the difference $\psi_{2}-\psi_{1}$ must be exact, and we use the method of \cite[Section ?]{McGowan} to build a globally defined form $\psi$ with $d\psi=\phi$ and with control over the Rayleigh-Ritz quotient.  Since our open cover has only two pieces, this produces no difficulties.  If we choose $y^{2}$ above in such a way that summing the relevant Rayleigh-Ritz inequalities gives the correct lower eigenvalue bound, we have
\begin{theorem} Suppose $N_\epsilon=\tilde{N}\cup (M^n\times I)$, $\tilde{N}$ and $M^n$ compact, with metric \begin{equation}g_\epsilon=\rho(\epsilon,t)^{2a}dt^2 +\rho(\epsilon,t)^{2b}ds_M^2\end{equation} on $M^n \times I$, with $\rho$ as above, and $d$ the dimension of the space of harmonic $p$-forms on $M^{n}$.  Let $$R=\int_0^1\rho(\epsilon,s)^a\,ds.$$  Let $N_\epsilon(x^2)$ be the number of eigenvalues of the Laplacian acting on coexact $p$-forms in $[\sigma,\sigma + x^2)$ where $\sigma$ is the bottom of the essential spectrum for coexact forms of degree $p$ and $0 < p < n$.  Then $$N_\epsilon(x^2) =\frac{dxR}{\pi} + O_x(1)$$ where $d$ is as in theorem \ref{th1}.\end{theorem}

In the special case when $d=0$, i.e. when there are no harmonic forms on the cross section, we have the following corollary,
\begin{corollary} Suppose $N$ is as above, with $d=0$ for some $p$.  Then the essential spectrum of the Laplacian acting on exact forms of degree $p$ is empty.\end{corollary}

 \end{document}